# Geographic Characterisation of Children's Health and Wellbeing through Vulnerability Indices Using Principal Component Analysis (PCA) and K-Means Clustering

Harriette Phillips[1*], Aiden Price[1], Kerrie Mengersen[1]

[1]Queensland University of Technology – School of Mathematical Sciences 2 George St

Brisbane City QLD 4000, Australia

*Corresponding author

Using Australian children as a case study, this paper aims to develop vulnerability indices to characterise the health and well-being of children aged 0-5 years old. The indices are used to identify differences in children's health and well-being across geographic regions and identify clusters of regions with similar characteristics. The approach is underpinned by two well-known statistical methods, namely Principal Component Analysis and K-Means clustering. The identification of these regions with similar vulnerability characteristics can then be used to derive new insights into drivers of children's health and well-being and support improved decision-making for services in Australia.

# 1. Introduction

Informative data assets are an essential foundation for evidence-based approaches to improving health and wellbeing for a given population (1,2). These data assets can facilitate targeted responses, better justify health and budget decisions, and improve transparency in policymaking. An effective way to develop evidence-based public health policies and decisions is through understanding the geographical differences in the health and wellbeing of communities (2,3). These differences can be captured through mapping, monitoring, and modelling a range of health, demographic, social and geographic data at a small area spatial scale.

One of the most important populations that require this assessment is children (4–6). The influence of where children live can have an impact on their health and wellbeing. Child health is also an important factor in the overall health of a population (4–6). Differences in regions can influence health vulnerability, where regions with better access to services and resources are likely to have better outcomes. Children who have not had their needs met are more vulnerable to harm (6). Vulnerability can be defined as children who are more susceptible to harm and/or receive less assistance (7–11). There is no single one cause for vulnerability: it is associated with a range of inequalities and complex interactions of economic, social, political, technological, and factors within communities. Since vulnerability cannot be directly observed, it is typically measured as a latent variable through one or more vulnerability indices. These indices provide powerful data assets that can aid in community level policy planning and the mitigation of current health challenges as well as reducing any future risks (12).

Literature on vulnerability indices for child health and wellbeing is relatively sparse. An early paper, co-developed with UNICEF and published in 2013, created a vulnerability index for adolescent girls in Uganda [16]. The study used a cut-off technique in which a girl was declared to be highly vulnerable if she fell below given thresholds at all three of the individual, family, and community dimensions. These dimensions were informed by two well-known national surveys. The percentage of highly vulnerable girls in each region was reported to be useful for more informed resource allocation. A more recent paper, published in 2023, reported on the relationship between a social vulnerability index and surgical outcomes for children in California, USA (13). The index was constructed by the Centre for Disease Control and represent census tracts across four main themes: housing and transportation, minority status and language, household composition and disability, and socioeconomic status. The vulnerability indices were rated using percentile rankings on a national scale, where higher percentiles indicated higher vulnerability. High vulnerability was defined when regions were above the 90th percentile. Australia, the focal country for the present paper, has been the subject of one relevant vulnerability index study. A recent paper, published in 2023, implemented k-means clustering to identify groups of vulnerability and the characteristics of these groups (14). The index used data from the Australian Early Development Census and the Australian Bureau of Statistics for Queensland, Australia. While the study highlights the link between vulnerability and geographic location this is limited to Queensland.

Various methodology can be used to develop vulnerability indices. Two exemplar methods were developed in the context of estimating flood vulnerability (15,16). The earlier study, published in 2016, uses Principal Component Analysis and K-Means clustering to develop vulnerability indices (15). The indices identify groups of vulnerable communities and profile the demographic characteristics of each group. The vulnerability indices were also categorised using z-scores into five different levels of vulnerability, from very low to very high. The more recent study, published in 2023, compares a weighted approach and a principal components approach to creating vulnerability indices, with the aim of prioritising actions before flood disasters occur. (16).

The current paper aims to develop a more comprehensive vulnerability index for children's health and wellbeing in Australia. The methodology builds on and extends previous approaches within and outside the field of child health [13], [17], [18]. First, principal component analysis (PCA) is implemented to create the indices as a weighted linear combination of the primary measured variables (16). K-Means clustering is then implemented to group regions across Australia, in order to characterise regions with common vulnerability characteristics, and identify spatial trends in the indices (17). The child health and wellbeing vulnerability indices are developed in the context of a larger endeavour to create a suite of publicly available data assets for the Australian Child and Youth Wellbeing Atlas (ACYWA) (18). This atlas has garnered support from over twenty-five industry partners including universities, and government agencies. The Atlas includes over 960 individual indicators collected across six domains: valued, loved, and safe; material basics; healthy; learning; participation; and positive sense of identity and culture (18,19). The Atlas aims to support better evidence-based policy, action and decision-making based on the latent geographical differences in communities, and in turn reduce the cost on health care, improve child health and wellbeing, and advance overall community health.

## 2. Methods

The development of a vulnerability index is comprised of two key components: data and methodology (20). The data component involves the identification of candidate variables, acquisition of corresponding data from relevant datasets, final selection of variables for inclusion in the index, and determination of the spatial resolution for index development and communication. The order of these steps may vary depending on the context. The methodology component involves statistical modelling of the data to create the indices, and spatial clustering of regions based on the derived indices. These components are described in detail below in the context of the case study described above.

# 3. Candidate Variables

This study focused on Australian children aged 0-5 years old during the period 2006-2022. Data were aggregated over period (2006-2022) and geographic region (at the Australian Statistical Geographical Standard at Statistical Area Level 3 (SA3) (21)). SA3 contains 358 non-intersecting regions across Australia (see Geographic Resolution and Appendix 1). The majority of the data were collected as part of the ACYWA study (18). Additional variables were collected from other free and publicly available resources including the Australian Early Development Census, the Australian Bureau of Statistics, and the Australian Curriculum, Assessment and Reporting Authority (22–24). Data collected were de-identified and aggregated, and participant information was not collected. Consequently, researchers did not have access to participant information at any time throughout the creation of this paper. These data were accessed from 28th of February 2023.

# 4. Variable Selection

A total of 41 variables were considered as candidates for inclusion in the indices. A list of the candidate variables, including a brief description, source and original creator, is provided in Appendix 2.

A variable selection process was then undertaken to reduce this set to a smaller number of most appropriate variables [18]. The final variables selected for use is outlined in Section 3.1. The variable selection process comprised the following steps.

Step 1: Regions that contained zero and NA for the estimated regional population (ERP) (25) were omitted from any further statistical methods and processes. Omitted regions are listed in Appendix 1. Variables with a high number of zero counts and not applicable (NA) entries for all variables were also considered for removal. For variables with a relatively small number of missing entries, the NA values were replaced with the mean of neighbouring regions. Regions were considered neighbours if they shared one or more shared boundary points (26). The `pol2ng` function from the `spdep` package was then applied in RStudio to find region neighbours based on their geographical boundaries (26).

Step 2: Normality of the input variables is a typical requirement, or at least assumed, for the component PCA and K-Means components of the vulnerability index development (15). Normality was assessed by measuring the skewness of the data through visualisations (boxplot, histogram, Q-Q plot) [50] and statistical thresholds: <|2| (acceptable) and >||2| (unacceptable) (27,28). Skewness values were calculated using `skewness` function from the `moments` package in RStudio (29). Variables with unacceptable skewness were all positively skewed so a log transformation was applied, and the transformed variables were assessed. If they were still of unacceptable skewness (exceeding ±2), they were considered for removal.

Step 3: Correlations between the variables that passed Steps 1 and 2 were then evaluated (30,31). Based on previous literature, variables with a very high Pearson correlation (± .90) were removed for use in the vulnerability indices (32,33).

# 5. Geographical Resolution

The Australian Statistical Geographical Standard (ASGS) was selected as the foundation for the spatial resolution of the vulnerability indices. The ASGS is a hierarchy of structures used to classify Australia that was first introduced in 2011 (21,34). These structures are updated every five years to reflect changes and growth in the Australian population and is used by many organisations to classify, collect, release, and analyse geographical data in a standardised format. This paper focuses on three of the structures outlined by the ASGS: Statistical Area 3 (SA3), Remoteness, and States and Territories. SA3 is the most granular of the structures investigated (35). These regions typically contain between 30 000 to 130 000 people. The SA3s reflect areas with similar social and economic characteristics. All variables were sourced at SA3. Remoteness areas are domains provided by the ABS and capture a regions access to services (36). This structure can be divided into five categories: Major Cities, Inner Regional, Outer Regional, Remote, and Very Remote. In this study, Statistical Area Level 1 remoteness area data were aggregated to SA3 to enable comparisons with derived vulnerability indices and their characteristics. This allows for trends and connections across regions to be investigated. Finally, the states were also investigated. Australia is broken up into nine states by the ASGS: New South Wales (NSW), Victoria (VIC), Queensland (QLD), South Australia (SA), Western Australia (WA), Tasmania (TAS), Northern Territory (NT), Australian Capital Territory (ACT), and Other Territories (OTHER) (i.e., Territory of the Cocos (Keeling) Islands and Norfolk Islands, Territory of Christmas Islands, Jervis Bay Territory) (21). Similarly to remoteness categories, the results of the vulnerability indices were compared across the different states.

An overview of the total number of regions in each state across the five domains of remoteness is provided in Table 1. Regions with no populations according to the ERP (n=5) and non-spatial regions (n=18) were excluded. The three spatial resolutions are visualised in Figure 1.

Table 1. The total number of regions per each remoteness domains from the Australian Statistical Geographical Standard (21).

| State | Major Cities | Inner Regional | Outer Regional | Remote | Very Remote | Total |
|-------|-------|--------|--------|--------|--------|-------|
| NSW | 52 | 25 | 10 | 2 | 1 | 90 |
| QLD | 49 | 17 | 12 | 2 | 1 | 81 |
| VIC | 40 | 20 | 6 | 0 | 0 | 66 |
| WA | 21 | 2 | 4 | 4 | 3 | 34 |
| SA | 18 | 3 | 5 | 1 | 0 | 27 |
| TAS | 0 | 7 | 7 | 1 | 0 | 15 |
| ACT | 9 | 1 | 0 | 0 | 0 | 10 |
| NT | 0 | 0 | 4 | 2 | 2 | 8 |
| OTHER | 0 | 1 | 0 | 0 | 3 | 4 |
| Total | 189 | 76 | 48 | 12 | 10 | 335 |

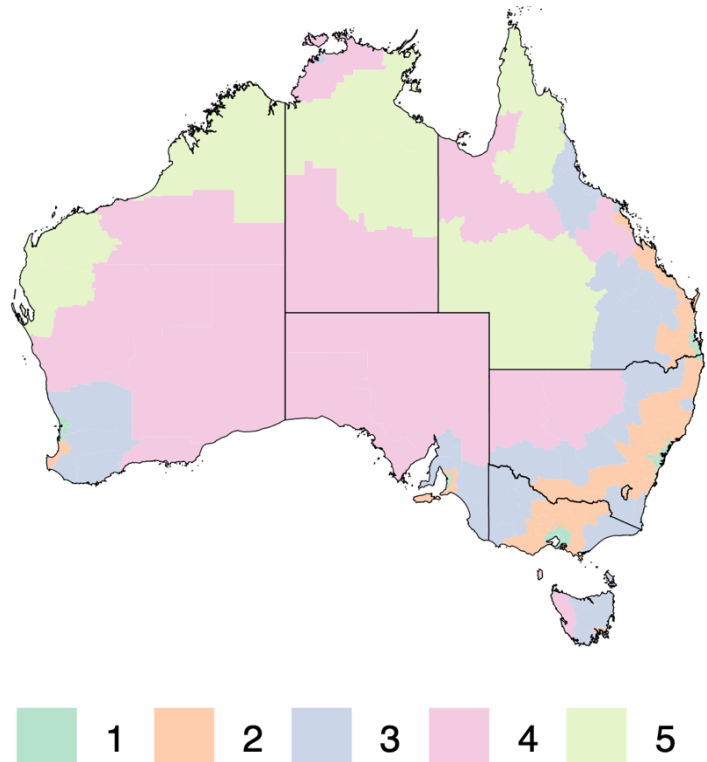

Figure 1. Visualisation of the remoteness domains across the Statistical Area Level Three regions of Australia (21). The legend refers to the domains of remoteness which are colour-coded with the corresponding values: 1: Major Cities, 2: Inner Regions, 3: Outer Regions, 4: Remote, and 5: Very Remote. The boundaries of each state are outlined in black.

# 6. Creation of the Vulnerability Indices

After the variables were selected, Principal Component Analysis (PCA) was implemented in RStudio to create the initial vulnerability indices [43]. PCA is a well-known dimension reduction technique that has become a commonly used approach to develop vulnerability indices (12,16,39,40). The outputs of the PCA are a set of eigenvectors and associated eigenvalues. The eigenvector describes the direction of data onto a reduced set of orthogonal axes, and the eigenvalues describe how the original data are transformed onto these new axes (39). Each eigenvector is a weighted sum of the original variables. The relative magnitudes of the weights were used to inform which variables contribute highly to each component (39).

The eigenvalues, which also inform about the proportion of total variance explained by each principal component, were used as a measure of its relative value or merit as a candidate vulnerability index.

Three approaches are commonly used to derive indices from principal components (16,41–44). The first is to retain the first component only as it explains the largest amount of variance within data (16,45). The second is to set a threshold for the cumulative variance and retain the number of components required to meet this threshold. Previous literature has used 30%, 50-80% as thresholds (41–44). The third approach is to keep components with eigenvalues close to or higher than one (12,15,16). This third approach was adopted for the current study. The `prcomp` function from the `stats` package in RStudio was used to create the PCA outputs (46). The retained principal components were declared to be the new derived vulnerability indices. This technique has previously been implemented in vulnerability index literature (12,15,16).

# 7. Clustering of Vulnerability Indices

The vulnerability indices obtained through the PCA step described in the previous subsection were obtained for each region of Australia. To understand the spatial profile of child health and wellbeing and identify spatial patterns, the regions were clustered with respect to their vulnerability indices. This was implemented in RStudio using the `kmeans` function from the `stats` package (17,47). K-Means clustering has been employed previously for vulnerability assessment in other fields such as flooding (14,48–50).

The K-Means analysis creates groups that are similar to each other with respect to the values of the input variables (in this case the vulnerability indices) whilst maximising the differences between the groups (15). K-Means clustering is achieved through the following steps (17,48,51,52). First, the number of clusters is defined by the user. The number of clusters is often chosen through statistical or visual methods (e.g., elbow method, information criterion approach, information theoretic approach, etc). Given an initial (often arbitrary) set of cluster centroids, regions are allocated to clusters based on the Euclidean distance to the cluster centroid. Centroids are then updated based on cluster formation. The process is repeated until the centroids do not change significantly.

# 8. Results

## 8.1 Variables

After the variable selection process, 26 variables were chosen for inclusion in the vulnerability index. These are described in Table 2 where a shortened version of the variable name is also included for later reference. The correlation between these variables is shown in Figure 2. All variables are continuous.

Table 2. Overview of variables used for the development of the vulnerability indices, shortened variable name, and brief description for each variable

| Variable | Short-Form | Definition |
|---|---|---|
| Estimated Resident Population | ERP | The estimate resident population provides information at a granular level about Australia's population, population growth, and population changes [25]. |
| Primary School | PRIM | The location of Australian primary schools [23] |
| Childcare Centres | CC | The location of childcare centres in Australia [53] |

| Healthcare Centres by Persons | HCC2 | Healthcare centre locations across Australia calculated using the population for a region. These healthcare centres include relevant services for children aged 0-5 [54]. Examples include Paediatric medical services, Paediatric dentistry services, and family counselling services [55]. |
| Developmentally Vulnerable on the Physical Health and Wellbeing Domain | DVH | A child's overall physical health and wellbeing is measured by their readiness for the school day, gross and fine motor skills as well as their physical independence [56]. |
| Developmentally Vulnerable on the Emotional Maturity Domain | DVE | Social competence for children is determined through several factors. These include a child's overall social competence, their approach to learning, responsibility, readiness to explore, and respect [56]. |
| Developmentally Vulnerable on the Social Competence Domain | DVS | A child's social competence is determined by their ability to learn, readiness to learn, as well as responsibility and respect of others [56]. |
| Developmentally Vulnerable on the Language and Cognitive Skills (School-Based) Domain | DVL | Children's literacy and numeracy, and interest in literacy is measured through this domain [56]. |

| Developmentally Vulnerable on the Communication Skills and General Knowledge Domain | DVC | A child's communication skills and general knowledge is measured in this domain [56]. |
|---|---|---|
| Developmentally Vulnerable on the Physical Readiness and Cleanliness for School Day Domain | DVR | Children were determined to be developmentally vulnerable based on their physical readiness and cleanliness of a child for the school day [56]. |
| Developmentally Vulnerable on the Physical Independence Domain | DVI | A child's physical independence is measured in this domain [56]. |
| Developmentally Vulnerable on the Gross and Fine Motor Skills Domain | DVG | A child's gross and fine motor skills is used to determine their developmental vulnerability in this domain[56]. |
| Small for Gestational Age | SML | The birthweight is less than 2750g when a baby is born on or after 40 weeks is considered small for gestational age [57], [58]. This is a key indicator for infant health used to measure both the health and wellbeing of mother and is a determinant of the baby's health, survival, development and health and wellbeing [58]. |
| Mother's Smoking During Pregnancy | SMO | Women who reported smoking during pregnancy. This is a preventable risk factor that results in adverse health outcomes for mothers and their babies [59]. |
| Five Plus Antenatal Visits | ANT | Antenatal visits during pregnancy exceeds five [57]. |
| First Antenatal Visit is Less than Fourteen Weeks | A14 | Woman whose first antenatal visit was during their first trimester (less than fourteen weeks of gestational age). Visits during the first trimester can be associated with improved maternal health during pregnancy, and health outcomes for the child [60]. |
| First Antenatal Visit is between fourteen and nineteen weeks | A19 | First antenatal visit is after the first trimester of pregnancy and is between fourteen and nineteen weeks [60]. |
| First Antenatal Visit is over Twenty weeks | A20 | Women whose first antenatal visit is over twenty weeks [60]. |
| Birth is less than thirty-seven weeks | BRTH | Babies who were born at less than thirty-seven weeks of gestations. Gestational age can be associated with a baby's overall health and generally poorer health outcomes are reported with babies that are born earlier [61]. |
| Mother's Age when giving birth is 15-19 | M15 | The age of the mother at time of birth is between 15-19. If the mother is under 20 this can be associated with adverse outcomes and increased risk of complications associated with pregnancy [62]. |
| Mother's Age when giving birth is 20-24 | M20 | The age of the mother at the time of birth is between 20-24 [62]. |
| Children who are Fully Immunised | IMM | Children who have received their immunised at 1year, 2 years, and 5 years old |
| Fertility Rate | FRT | The number of babies born per number of woman [63] |
| Mortality Rate per 100000 | MRT2 | The deaths of children aged 0 to 9 per 100000 [64]. |
| Preschool Attendance | ATT | Children that are aged four or five that have been enrolled in a preschool program [65]. |
| Core Activity Need for Assistance | CORE | Children who have a severe or profound core activity limitation and thus need assistance in their day-to-day lives. This must limit three aspects of their daily activities (e.g., self-care, communication, or mobility) and can be long term disabilities or long-term health conditions [66]. |

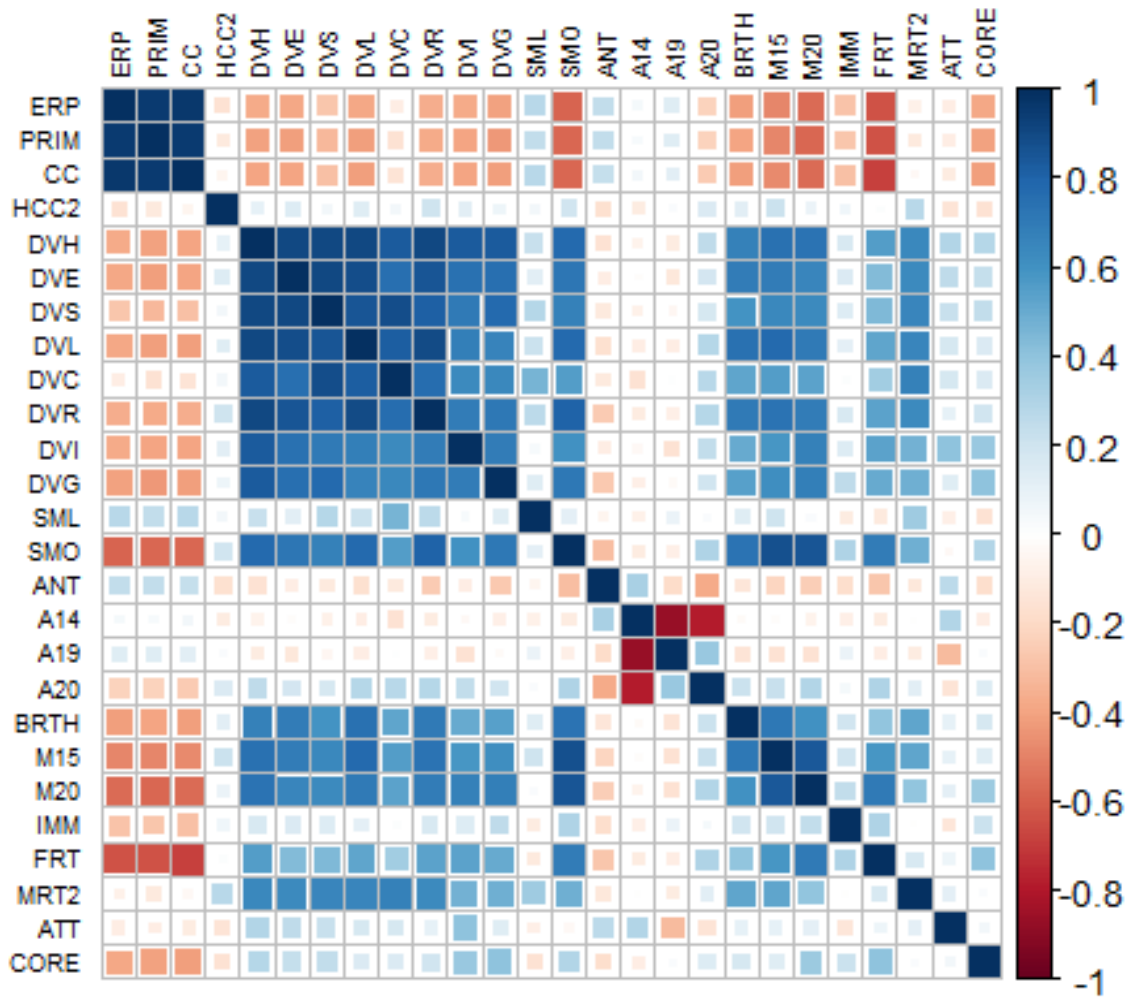

Figure 2. Pearson correlation plot for the variables selected for use in the development of the vulnerability indices. Strong positive correlations are represented by dark blue, strong negative correlation by dark red, and weak correlation by light to white colours. The list of variable names is available in Table 2.

Figure 2 highlights a few key correlations between the variables. ERP, PRIM, and CC are highly positively correlated with one another and have varying degrees of negative correlations with the majority of the remaining variables. SML, ANT, and A19 have weak positive correlations with ERP, PRIM, and CC. Additionally, A14 has weak to no correlation with ERP, PRIM, and CC. DVS, DVL, DVC, DVR, DVI, and DVG are all positively correlated with one another. A14, A19, and A20 are all negatively correlated with one another.

## 8.2 PCA Results

The principal component analysis identified five PCs with eigenvalues greater than one. These components explain 75.94% of the cumulative variance within the data. The proportion of variance, cumulative variance and eigenvalues for all components are provided in Appendix 3. Variables with a weight greater than ±0.20 were considered to contribute substantively to vulnerability. Five new vulnerability indices were constructed based on these results (Table 3). The statistical distribution of these indices was examined to confirm acceptable levels of skewness, potential outliers in the original data, and interpretability. The results of the skewness assessment are presented in Appendix 4.

Table 3. The first five vulnerability indices with eigenvalues greater than one. The eigenvalues and proportion of variance for these indices have been supplied in this table

|  | VI1 | VI2 | VI3 | VI4 | VI5 |
|---|---|---|---|---|---|
| Percent of Variance | 43.05 | 12.62 | 10.52 | 5.50 | 4.24 |
| Eigenvalues | 11.19 | 3.28 | 2.74 | 1.43 | 1.10 |

The first vulnerability index (VI1) explained the largest proportion of variance (43.05%) of the original data and is defined by the positive weights of two groups of variables: education and pregnancy (see Appendix 3). These include developmentally vulnerable education variables across physical health and wellbeing (DVH), language and cognitive skills (DVL), physical readiness and cleanliness for school day (DVR), emotional maturity (DVE), social competence (DVS), gross and fine motor skills (DVG), physical independence (DVI), and communication skills and general knowledge (DVC). The pregnancy variables include smoking during pregnancy (SMO), mother's age (20-24) (M20), mother's age (15-19) (M15), child born less than thirty-seven weeks (BRTH), and fertility rate (FRT).

The second vulnerability index (VI2) explained 12.62% of the variance in the original data and is positively weighted towards density and general health and wellbeing variables (see Appendix 3). There are three density dimensions: childcare centres (CC), Estimated Regional Population (ERP), and primary schools (PRIM). The remaining health and wellbeing variables include small for gestational age (SML), mortality rate (MRT2), communication skills and general knowledge (DVC), fertility rate (FRT), and need for assistance (CORE).

The remaining vulnerability indices had a combined variance of 20.26%. The dimensions of these indices include pregnancy variables and preschool attendance for VI3, pregnancy variables and general health and wellbeing variables for VI4 and VI5 (see Appendix 3). The non-standardised visualisation of the vulnerability indices has been provided in Figure 3. The individual vulnerability indices have different extremes as such, an additional visualisation of the indices on a consistent truncated scale has been provided in Appendix 3.

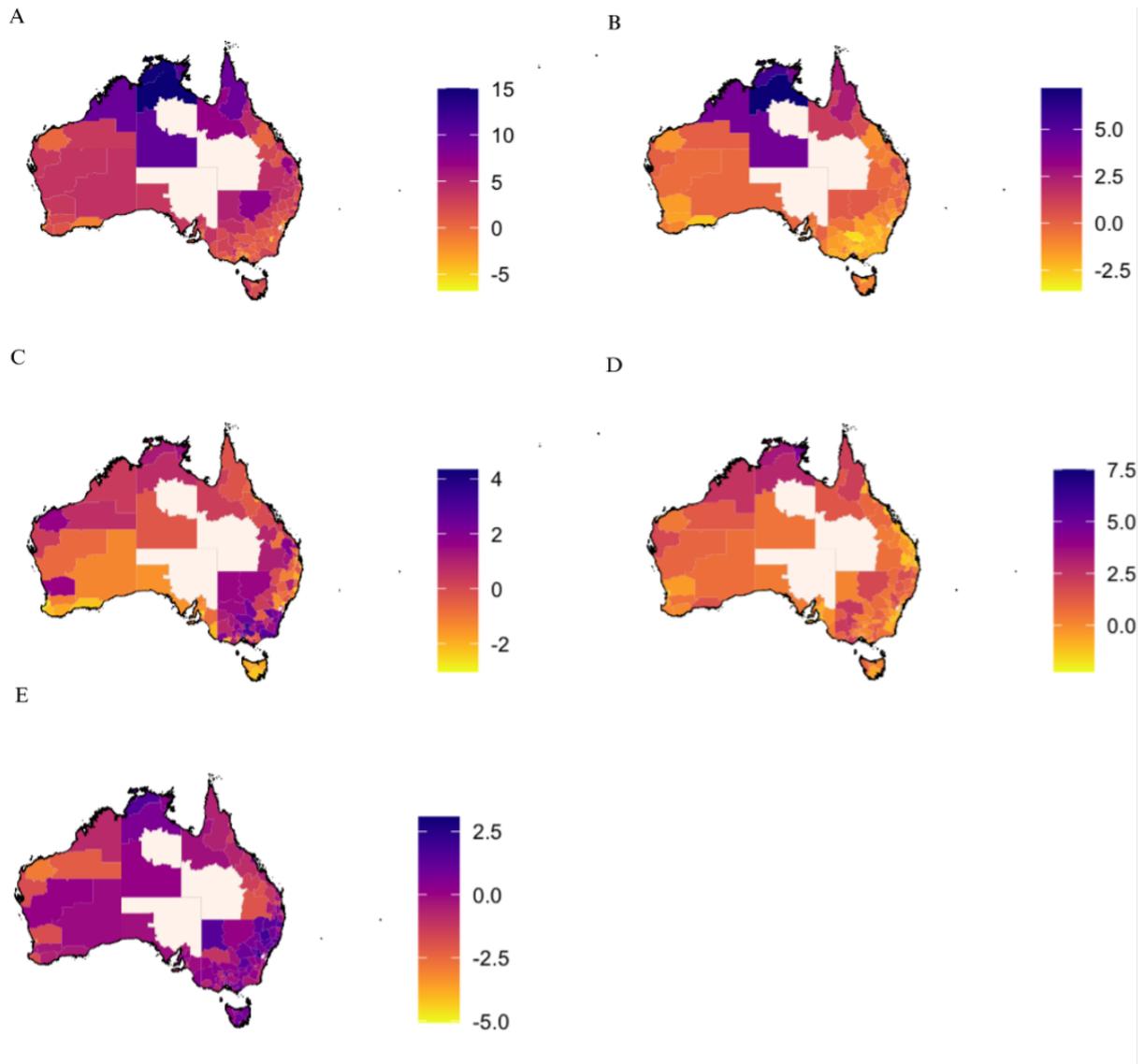

Figure 3 The alphabetical order in Figure 3 corresponds to the first five vulnerability indices in numerical order. These maps showcase each index on a non-standardised scale. Regions with no population based on the ERP are coloured white (see Appendix 1). Purple (high positive values) indicate high levels of vulnerability and yellow (high negative values) represent levels of low vulnerability.

## 8.3  K-Means Results

Following the implementation of PCA, K-Means clustering was applied to the five selected principal components (vulnerability indices). Four clusters were selected for the k-means implementation and were chosen by visual inspection of the scree plot using the elbow method (Appendix 5). The resulting cluster centroids and standard deviations are presented in Table 4 and a summary of the remoteness domains captured in each cluster is provided in Table 5. A box plot and choropleth map of the clusters are also presented in Figure 4 and Figure 5. The weights of the individual variables were also used to determine levels of vulnerability (Appendix 3).

Table 4. The vulnerability indices cluster centroids (mean) and standard deviations (sd).

|  | VI1 | VI2 | VI3 | VI4 | VI5 |
|---|---|---|---|---|---|
| Cluster | mean (sd) | mean (sd) | mean (sd) | mean (sd) | mean (sd) |
| C1 | 2.63(1.39) | -0.16(1.05) | -0.25(1.48) | -0.25(1.06) | 0.08(0.89) |
| C2 | 10.07(3.01) | 3.34(2.30) | 0.75(0.79) | 2.07(1.47) | -0.04(0.81) |
| C3 | -3.39(1.17) | 1.78(1.51) | 0.07(0.16) | 0.18(-0.04) | -0.11(0.00) |
| C4 | -1.01(1.59) | -1.46(0.91) | 0.16(1.81) | -0.04(1.13) | 0.00(1.11) |

Table 5. Distribution of vulnerability index clusters (C1 to C4) across the remoteness domains (number of regions per cluster).

|  | C1 | C2 | C3 | C4 |
|---|---|---|---|---|
| Major Cities | 37 | 0 | 86 | 66 |
| Inner Regional | 40 | 0 | 1 | 35 |
| Outer Regional | 39 | 1 | 1 | 7 |
| Remote | 6 | 4 | 0 | 2 |
| Very Remote | 1 | 4 | 0 | 5 |
| Total | 123 | 9 | 88 | 115 |

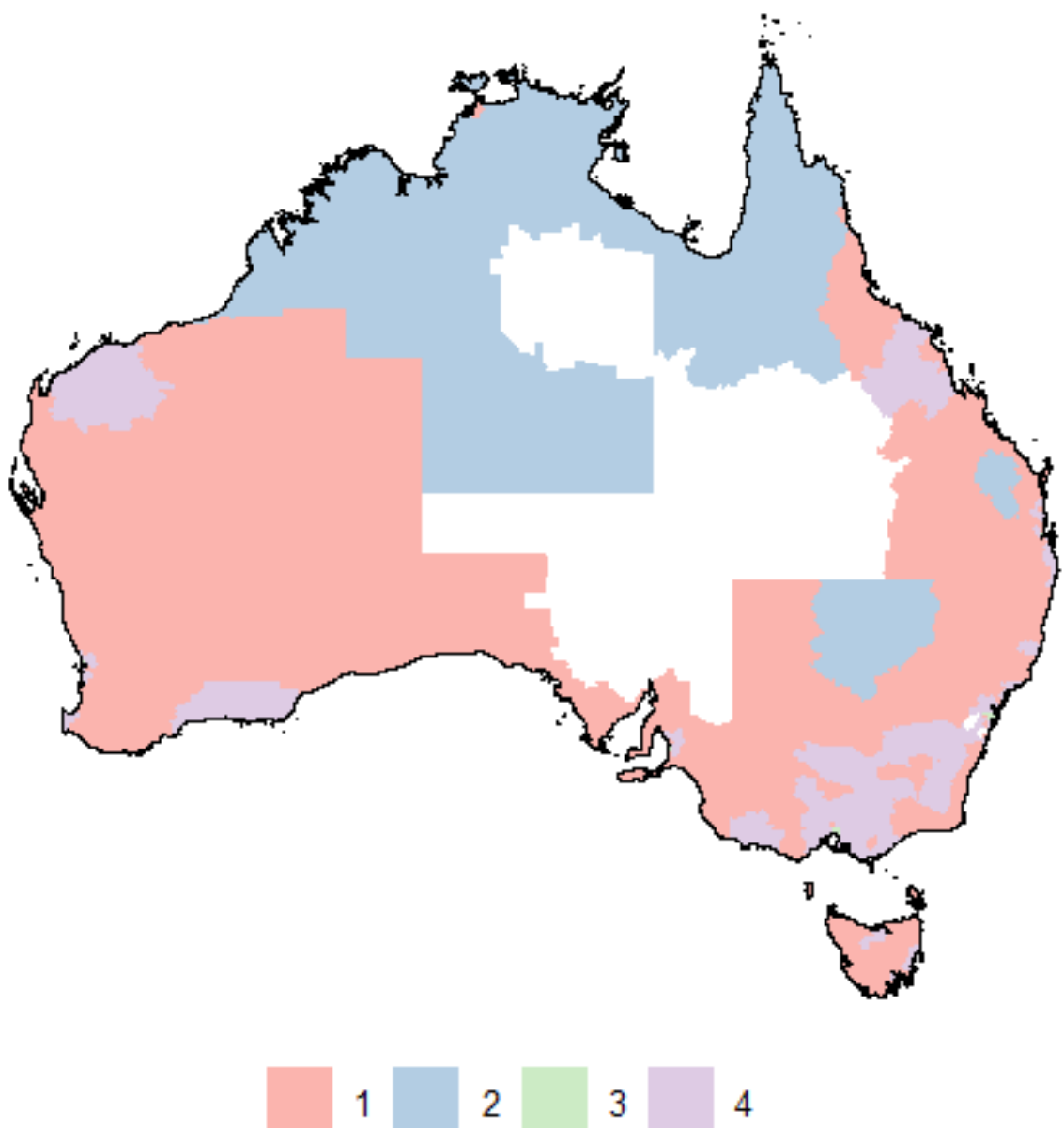

Figure 4. Visualisation of all clusters. Each cluster is colour coded for clarity. Regions with no children are not included in the colours corresponding to each cluster (see Appendix 1).

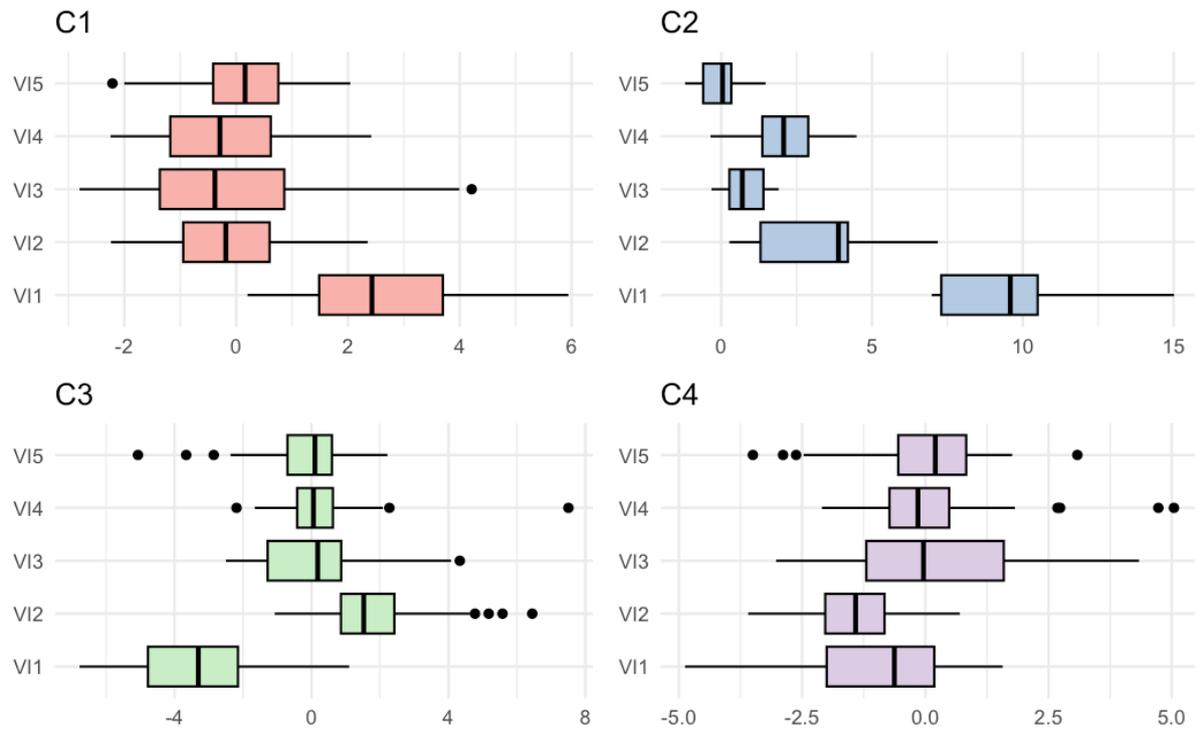

Figure 5. Boxplot visualisation the distribution of the vulnerability indices (VI1-VI5) of each cluster (C1-C4)

Within the clusters, four key groups of variables were identified: education, pregnancy, healthcare services, and general health and wellbeing. These variable grouping are used to describe the overarching themes within the clusters and subsequent vulnerability. Cluster 1 (C1) comprises 123 regions spanning all remoteness domains, with the majority in Major Cities, Inner Regional and Outer Regional regions (Table 5). This cluster is largely informed by education and pregnancy variables. Cluster 2 (C2) comprises 9 regions, mostly in Remote and Very Remote areas. All variable groups are present in this cluster. Cluster 3 (C3) comprises 88 regions, mostly in Major Cities. High vulnerability for C3 is mostly informed by education and pregnancy variables. Cluster 4 (C4) comprises 136 regions, mostly in Major Cities and Inner Regional areas. An overview of the variables that characterise the cluster profiles is presented in the Discussion. This cluster is also mostly informed by education and pregnancy variables. More specific details of these variables have been outlined in the Discussion.

# 9. Discussion

This paper has identified trends and characteristics of vulnerability by implementing PCA and K-Means clustering. Five distinct vulnerability indices were obtained, as well as four spatial clusters of regions across Australia based on these indices. These indices can help to provide insights into key factors associated of vulnerability; and highlight common trends, characteristics, and geographical differences across regions (2,3,53). For Australian children aged 0-5, several key trends and characteristics of vulnerability have been identified.

Vulnerability index 1 and Vulnerability index 2 can be used to capture vulnerability in all four clusters. These two indices provide the most insights into vulnerability across the clusters. Vulnerability Index 3 describes vulnerability in three clusters. Vulnerability Index 4 and Vulnerability Index 5 is used to capture vulnerability in two clusters.

Key insights derived from cluster characteristics provide regional insights, can be used to inform policy decisions, as well as provide targeted approached for an individual region's improving child health and wellbeing. For example, regions in Cluster 1 have good access to health care services (e.g., childcare centres, primary schools, and healthcare centres) and have their first antenatal visit between fourteen to nineteen weeks. These regions also have less babies that are small for their gestational age, and a have a lower mortality rate. There are two key areas of concern in this cluster: education and pregnancy variables. This includes a larger proportion of premature babies, mother's smoking during pregnancy, young mother's aged 15-24, as well as a lower fertility rate for regions in this cluster. There is also a larger proportion of children not performing well across education domains (e.g., language and cognitive skills, and communication skills and general knowledge). Preschool attendance is also less likely. Regions in this cluster are present across all remoteness domains with the majority of regions being Major Cities and their surrounding regions. If access to healthcare services declines, there could be potential risks to maternal and infant health, highlighting the need for improvements in both healthcare and education.

Further insights into vulnerability can be derived from the cluster characters for Cluster 2. Regions in this cluster also have two key areas of concern: education and pregnancy. These regions represent the worst performing regions across the education domains (e.g., language and cognitive skills, and communication skills and general knowledge). Furthermore, these regions have a significantly higher proportion of premature babies, mother's smoking during pregnancy, mother's aged 15-24, and children small for their gestational age, and the fertility rate is lower. Mothers in this cluster often have their first antenatal visit above 19 weeks and often have a minimum of five antenatal visits. Finally, the regions do not have access to healthcare services e.g., childcare centres, primary schools, and healthcare centres). The regions in this cluster are representative of remote regions in Australia. Overall, these regions face multiple vulnerabilities, especially in healthcare, education, and pregnancy, reflecting the challenges of remote regions.

A key area for concern for Cluster 3 is education (e.g., language and cognitive skills, and communication skills and general knowledge). The regions in this cluster also do not have good access to childcare centres and primary schools. Furthermore, there is a higher proportion of premature babies, mother's smoking during pregnancy, mother's aged 15-24, and children small for their gestational age as well as a lower fertility rate. The mortality rate is higher for regions in this cluster. The regions in this cluster are solely representative of Major Cities and two surrounding regions. Cluster 3 reflect the areas of need in urban areas.

Finally, Cluster 4 is performing well across all education domains (e.g., language and cognitive skills, and communication skills and general knowledge). These regions also have good access to primary schools and childcare centres. There is lower proportion of premature babies, mother's smoking during pregnancy, and young mothers aged 15-24, children who a small for their gestation age, and a higher fertility rate. Mother in these regions typically have their first antenatal first at 19 weeks or above and have at least five antenatal visits. These regions do not have good access for those who need core activity assistance. The regions in this cluster are present in all five remoteness domains. Although, most regions are Major Cities and their surrounding regions. Overall, this cluster is performing well and could build upon its current strengths to encouraged earlier antenatal visits and better core activity assistance. This would further enhance maternal and child health outcomes, building on their strengths in education and low rates of smoking and premature births.

There are a few limitations associated with the current research. First, to capture variations in vulnerability across the individual regions, decision makers may need to use their own judgment on if the current methodology is appropriate. This includes the selection of variables which can create different outputs depending on the variables selected. Additionally, this paper aims to develop vulnerability index for children's health and wellbeing which can inform implementation of new policies and decisions. It is important to note that health and wellbeing is not guaranteed where lowering vulnerability may not necessarily improve the health and wellbeing of all children even in optimal conditions (54). Additionally, there are a wide range of factors that contribute to the health and wellbeing of Australian children. These factors can be difficult to measure due to differences in procedures for collection and release and ethical and privacy impacts on sharing of many important data sources in child health. As a result, some variables are not available or publicly released, and those that are do not use a consistent format (e.g., measuring across different temporal, spatial, or demographics). As a result, even data that is available may not be usable for child health and wellbeing assessments. This limitation is being addressed through multiple organisations and initiatives, including ACYWA (18) by promoting improved data sharing practices and by working with data custodians to increase the volume free and publicly available data in the child health space. Although efforts exist to improve data availability, ongoing efforts are required from both government and research agencies to improve standardised availability of actionable data.

Another limitation of this approach is its sensitivity to changes to the starting seed when implementing K-Means. To assess the robustness of the results, four seeds (1767, 7462, 944, and 3401) were compared alongside the original seed (123). Using the Hubert-Arabie Adjusted Rand Index (ARI), it was found that many seed combinations (e.g., 123 versus 1767) were sensitive to changes. However, some pairs, such as 1767 and 3401, and 123 and 944 yielded more stable results. In previous literature, an ARI value below 0.65 indicated poor sensitivity (55). To accommodate this sensitivity, investigating methodology previously implemented in literature in literature could be explored. For example, a recent study explored K-Means ++ to compensate for sensitivity to changes to seeds (56).

The current work can be extended to include the comparison of the current methodology with other common statistical methodology used in vulnerability index development. This includes factor analysis, a statistical method which aims to explore the underlying patterns in data and simply complex relations between variables (57,58). Other implementations can include the comparing the two alternative approaches to principal component analysis (cut-offs and first component only) to the current approach (see Creation of Vulnerability Indices).

Overall, these vulnerability indices help address individual challenges for children by investigating latent geographical information (2,3,53). In turn, allowing for the needs of the overall population to be addressed (4–6). The complex interactions, measured through latent variables, have provided a powerful visualisation and profile of the characteristics of vulnerability for regions describing children's health and wellbeing in Australia that was not previously available. These indices support the work previously provided by the Australian

Child and Youth Wellbeing Atlas which aims to create evidence-based actions, policies and decision-making by providing free and publicly available child and youth health and wellbeing data (18). This knowledge can facilitate evidence-based justification and transparency in the implementation of new health and wellbeing strategies and policies. Evidence-based approaches can also help reduce the cost to health care systems and have a positive effect on children and future generations which benefits the whole of Australia (4,5,7).

# Acknowledgements

The current work was supported by the Centre for Data Science at the Queensland University of Technology, and a Research Training Program Stipend Scholarship awarded to Harriette Phillips. This work would also like to acknowledge Bernie Morrie, from the Australian Research Alliance for Children and Youth, for her guidance with the aims and structure in the development of the vulnerability indices. The Australian Child and Youth Wellbeing Atlas has also provided their support in the supply of the data for this project.

# Author Contributions

Ms Harriette Phillips: Conceptualisation, methodology, writing original draft, investigation, validation, project administration, software, data curation, writing review, and editing

Dr Aiden Price: Support in conceptualisation, methodology, writing review and editing, data curation, and supervision

Distinguished Professor Kerrie Mengersen: Support in conceptualisation, methodology, writing review and editing, and supervision

# 11. Appendix

Appendix 1. Regions omitted based on zero or not applicable entries for estimated regional population

Five regions were omitted from use in further statistical methodology based on Estimated Regional Population (ERP) zero entries as shown in Table 5. There are 2 from New South Wales, 1 from Queensland, 1 from South Australia, and 1 from the Northern Territory. There were an additional 18 regions remoted based on NA counts in ERP and these regions being classified as non-spatial regions (see Geographical Resolution) (35,59). These regions cannot be visualised regardless of their entries as non-spatial regions do not contain the necessary geographical boundaries.

Table 6. Regions and corresponding states for the zero counts and not applicable counts based on estimated regional population

| Zero Counts | |
|---|---|
| Statistical Area Three Region | State Name |
| Illawarra Catchment Reserve | New South Wales |
| Blue Mountains - South | New South Wales |
| Outback - South | Queensland |
| Outback - North and East | South Australia |
| Barkly | Northern Territory |
| Not Applicable Counts | |
| Migratory - Offshore - Shipping (NSW) | New South Wales |
| No usual address (NSW) | New South Wales |
| Migratory - Offshore - Shipping (VIC) | Victoria |
| No usual address (VIC) | Victoria |
| Migratory - Offshore - Shipping (QLD) | Queensland |
| No usual address (QLD) | Queensland |
| Migratory - Offshore - Shipping (SA) | South Australia |
| No usual address (SA) | South Australia |
| Migratory - Offshore - Shipping (WA) | Western Australia |
| No usual address (WA) | Western Australia |
| Migratory - Offshore - Shipping (TAS) | Tasmania |
| No usual address (TAS) | Tasmania |
| Migratory - Offshore - Shipping (NT) | Northern Territory |
| No usual address (NT) | Northern Territory |
| Migratory - Offshore - Shipping (ACT) | Australian Capital Territory |
| No usual address (ACT) | Australian Capital Territory |
| Migratory - Offshore - Shipping (OT) | Other Territories |
| No usual address (OT) | Other Territories |

Appendix 2. An overview of the variables used in the development of vulnerability indices

| Variable | Short-Form | Definition | Sourced From | Creator of Dataset |
|---|---|---|---|---|
| Estimated Resident Population | ERP | The estimate resident population provides information at a granular level about Australia's population, population growth and change (25). | Australian Child Youth and Well-being Atlas (18) | Australian Bureau of Statistics (61) |
| Primary School by Area | PRIM | The location of Australian primary schools calculated using the area of a given region (23). | Australian Curriculum, Assessment and Reporting Authority (23) | Australian Curriculum, Assessment and Reporting Authority (62) |
| Primary School by Persons | PRIM2 | The location of Australian primary schools calculated using the ERP for a given region (23). | Australian Curriculum, Assessment and Reporting Authority (23) | Australian Curriculum, Assessment and Reporting Authority (62) |
| Childcare Centres by Area | CC | Calculated using the population for a given region is the location of childcare centres in Australia (63) | The Australian Children's Education and Care Quality Authority (ACECQA) (63) | The Australian Children's Education and Care Quality Authority (ACECQA) (63) |
| Childcare Centres by Persons | CC2 | Calculated using the area for a given region is the location of childcare centres in Australia (63) | The Australian Children's Education and Care Quality Authority (ACECQA) (63) | The Australian Children's Education and Care Quality Authority (ACECQA) (63) |
| Hospitals by Area | HOS | The location of private and public hospitals in Australia by area in a given region (64). | Australian Institute of Health and Welfare (64) | Australian Institute of Health and Welfare (64) |

| Hospitals by Persons | HOS2 | The location of private and public hospitals in Australia by population of a given region (64). | Australian Institute of Health and Welfare (64) | Australian Institute of Health and Welfare (64) |
|---|---|---|---|---|
| Healthcare Centres by Area | HCC | Healthcare centre locations across Australia calculated using the given area for a region. These healthcare centres include relevant services for children aged 0-5 (65). Examples include Paediatric medical services, Paediatric dentistry services, and family counselling services (66). | National Health Services Directory (66) | National Health Services Directory (66) |
| Healthcare Centres by Persons | HCC2 | Healthcare centre locations across Australia calculated using the population for a region. These healthcare centres include relevant services for children aged 0-5 (65). Examples include Paediatric medical services, Paediatric dentistry services, and family counselling services (66). | National Health Services Directory (66) | National Health Services Directory (66) |
| Developmentally Vulnerable on One or More Domains | DV1 | A summary indicator used to determine whether a child is developmentally vulnerability on one or more | Australian Early Development Census (22) | Australian Early Development Census (22) |

| | | domains: physical health and well-being, social competence, emotional maturity, school-based language and cognitive skills, or communication skills and general knowledge (67) | | |
|---|---|---|---|---|
| Developmentally Vulnerable on Two or More Domains | DV2 | A summary indicator used to determine whether a child is developmentally vulnerability on two or more domain: physical health and well-being, social competence, emotional maturity, school-based language and cognitive skills, or communication skills and general knowledge (67). | Australian Early Development Census (22) | Australian Early Development Census (22) |
| Developmentally Vulnerable on the Physical Health and Well-being Domain | DVH | A child's overall physical health and well-being is measured by their readiness for the school day, gross and fine motor skills as well as their physical independence (67). | Australian Early Development Census (22) | Australian Early Development Census (22) |
| Developmentally Vulnerable on the Emotional Maturity Domain | DVE | Social competence for children is determined through several factors. These include a child's | Australian Early Development Census (22) | Australian Early Development Census (22) |

| | | overall social competence, their approach to learning, responsibility, readiness to explore, and respect (67). | | |
|---|---|---|---|---|
| Developmentally Vulnerable on the Social Competence Domain | DVS | Developmentally vulnerable based on relevant question and cut-off for social competence, learning and readiness, as well as responsibility and respect (67). | Australian Early Development Census (22) | Australian Early Development Census (22) |
| Developmentally Vulnerable on the Language and Cognitive Skills (School-Based) Domain | DVL | Children's literacy and numeracy, interest in literacy is measured. This is then determined to be developmentally vulnerable based off of cut-offs and relevant questions (67). | Australian Early Development Census (22) | Australian Early Development Census (22) |
| Developmentally Vulnerable on the Communication Skills and General Knowledge Domain | DVC | Relevant questions and cut-offs were used to determine if a child's communication skills and general knowledge was developmentally vulnerable (67). | Australian Early Development Census (22) | Australian Early Development Census (22) |
| Developmentally Vulnerable on the Physical Readiness and Cleanliness for School Day Domain | DVR | The physical readiness and cleanliness of a child for the school day. Children were determined to be developmentally vulnerable based | Australian Early Development Census (22) | Australian Early Development Census (22) |

| | | off of cut-offs and relevant questions (67). | | |
|---|---|---|---|---|
| Developmentally Vulnerable on the Physical Independence Domain | DVI | Cut-offs and relevant questions were used to determine if child's physical independence was developmentally vulnerable (67). | Australian Early Development Census (22) | Australian Early Development Census (22) |
| Developmentally Vulnerable on the Gross and Fine Motor Skills Domain | DVG | A child was determined to be developmentally vulnerable for their gross and fine motor skills based on cut-off and relevant questions (67). | Australian Early Development Census (22) | Australian Early Development Census (22) |
| Live Low Birth Weights | LOW | Live births where the weight of the baby does not exceed 2500g (68). Birthweight can be used as a determinant for the immediate and future health and well-being of a baby (69). | National Perinatal Data Collection (70) | Australian Institute of Health and Welfare (National Perinatal Data Collection) (70) |
| Small for Gestational Age | SML | The birthweight is less than 2750g when a baby is born on or after 40 weeks is considered small for gestational age (68,71). This is a key indicator for infant health used to measure both the health and well-being of mother and is a determinant of the baby's health, survival, | National Perinatal Data Collection (70) | Australian Institute of Health and Welfare (National Perinatal Data Collection) (70) |

| | | development and health and well-being (71). | | |
|---|---|---|---|---|
| Mother's Smoking During Pregnancy | SMO | Women who reported smoking during pregnancy. This is a preventable risk factor that results in adverse health outcomes for mothers and their babies (72). | Australian Child Youth and Well-being Atlas (18) | Australian Institute of Health and Welfare (National Perinatal Data Collection) (70) |
| Five Plus Antenatal Visits | ANT | Antenatal visits during pregnancy exceeds five (68). | National Perinatal Data Collection (70) | Australian Institute of Health and Welfare (National Perinatal Data Collection) (70) |
| First Antenatal Visit is Less than Fourteen Weeks | A14 | Woman whose first antenatal visit was during their first trimester (less than fourteen weeks of gestational age). Visits during the first trimester can be associated with improved maternal health during pregnancy, and health outcomes for the child (73). | Australian Child Youth and Well-being Atlas (18) | Australian Institute of Health and Welfare (National Perinatal Data Collection) (70) |
| First Antenatal Visit is between fourteen and nineteen weeks | A19 | First antenatal visit is after the first trimester of pregnancy and is between fourteen and nineteen weeks (73). | Australian Child Youth and Well-being Atlas (18) | Australian Institute of Health and Welfare (National Perinatal Data Collection) (70) |
| First Antenatal Visit is over Twenty weeks | A20 | Women whose first antenatal visit is over twenty weeks (73). | Australian Child Youth and Well-being Atlas (18) | Australian Institute of Health and Welfare |

| | | | | (National Perinatal Data Collection) (70) |
|---|---|---|---|---|
| Birth less than thirty-seven weeks | BRTH | Babies who were born at less than thirty-seven weeks of gestations. Gestational age can be associated with a baby's overall health and generally poorer health outcomes are reported with babies that are born earlier (74). | Australian Child Youth and Well-being Atlas (18) | Australian Institute of Health and Welfare (National Perinatal Data Collection) (70) |
| Mother's Age when giving birth is 15-19 | M15 | The age of the mother at time of birth is between 15-19. If the mother is under 20 this can be associated with adverse outcomes and increased risk of complications associated with pregnancy (75). | Australian Child Youth and Well-being Atlas (18) | Australian Institute of Health and Welfare (National Perinatal Data Collection) (70) |
| Mother's Age when giving birth is 20-24 | M20 | The age of the mother at the time of birth is between 20-24 (75). | Australian Child Youth and Well-being Atlas (18) | Australian Institute of Health and Welfare (National Perinatal Data Collection) (70) |
| Children who are Fully Immunised | IMM | Children who have received their immunised at 1year, 2 years, and 5 years old | Australian Child Youth and Well-being Atlas (18) | Australian Government Department of Health and Aged Care (76) |
| Fertility Rate | FRT | The number of babies born per number of woman (77) | Australian Child Youth and Well-being Atlas (18) | Australian Bureau of Statistics (77) |

| Mortality Rate per 1000 | MRT | The deaths of infants (under the age of one) per 1000 (78). | Australian Child Youth and Well-being Atlas (18) | Australian Institute of Health and Welfare (National Mortality Database) (78) |
|---|---|---|---|---|
| Mortality Rate per 100000 | MRT2 | The deaths of children aged 0 to 9 per 100000 (78). | Australian Child Youth and Well-being Atlas (18) | Australian Institute of Health and Welfare (National Mortality Database) (78) |
| Preschool Attendance | ATT | Children that are aged four or five that have been enrolled in a preschool program (79). | Australian Child Youth and Well-being Atlas (18) | Australian Bureau of Statistics (Preschool Education) (79) |
| Core Activity Need for Assistance | CORE | Children who have a severe or profound core activity limitation and thus need assistance in their day-to-day lives. This must limit three aspects of their daily activities (e.g., self-care, communication, or mobility) and can be long term disabilities or long-term health conditions (80). | Australian Child Youth and Well-being Atlas (18) | Australian Bureau of Statistics (Census) (81) |
| Total Homelessness | HMLS | Homeless children living in non-convention, emergency, or short-term accommodation (82). Australian children experiencing homelessness are | Australian Institute of Health and Welfare (Specialist Homelessness Services Collection) (83) | Australian Institute of Health and Welfare (64) |
| | | amongst the most socially and economically disadvantaged. Australian Child Youth and Well-being Atlas (18) | | |
| Homelessness As a Result of Accommodation Reasons | ACC | The reported reason for homelessness and seeking assistance is accommodation (83). | Australian Child Youth and Well-being Atlas (18) | Australian Institute of Health and Welfare (Specialist Homelessness Services Collection) (83) |
| Homelessness As a Result of Financial Reasons | FIN | The reported reason for homelessness and seeking assistance is financial (83). | Australian Child Youth and Well-being Atlas (18) | Australian Institute of Health and Welfare (Specialist Homelessness Services Collection) (83) |
| Homelessness As a Result of Health Reasons | HLTH | The reported reason for homelessness and seeking assistance is health (83). | Australian Child Youth and Well-being Atlas (18) | Australian Institute of Health and Welfare (Specialist Homelessness Services Collection) (83) |
| Homelessness As a Result of Interpersonal Relationships | IR | The reported reason for homelessness and seeking assistance is interpersonal relationships (83). | Australian Child Youth and Well-being Atlas (18) | Australian Institute of Health and Welfare (Specialist Homelessness Services Collection) (83) |
| Couples Who Are Homeless with One Or More Children | CPLS | Clients reporting homelessness are couples with one or more children (83). | Australian Child Youth and Well-being Atlas (18) | Australian Institute of Health and Welfare (Specialist Homelessness Services Collection) (83) |

| Variable | Short-Form | Definition | Sourced From | Creator of Dataset |
|---|---|---|---|---|
| Estimated Resident Population | ERP | The estimate resident population provides information at a granular level about Australia's population, population growth and change [1]. | Australian Child Youth and Well-being Atlas [2] | Australian Bureau of Statistics [3] |
| Primary School by Area | PRIM | The location of Australian primary schools calculated using the area of a given region [4]. | Australian Curriculum, Assessment and Reporting Authority [4] | Australian Curriculum, Assessment and Reporting Authority [5] |
| Primary School by Persons | PRIM2 | The location of Australian primary schools calculated using the ERP for a given region [4]. | Australian Curriculum, Assessment and Reporting Authority [4] | Australian Curriculum, Assessment and Reporting Authority [5] |
| Childcare Centres by Area | CC | Calculated using the population for a given region is the location of childcare centres in Australia [6] | The Australian Children's Education and Care Quality Authority (ACECQA) [6] | The Australian Children's Education and Care Quality Authority (ACECQA) [6] |
| Childcare Centres by Persons | CC2 | Calculated using the area for a given region is the location of childcare centres in Australia [6] | The Australian Children's Education and Care Quality Authority (ACECQA) [6] | The Australian Children's Education and Care Quality Authority (ACECQA) [6] |
| Hospitals by Area | HOS | The location of private and public hospitals in Australia by area in a given region [7]. | Australian Institute of Health and Welfare [7] | Australian Institute of Health and Welfare [7] |
| Hospitals by Persons | HOS2 | The location of private and public hospitals in Australia by | Australian Institute of Health and Welfare [7] | Australian Institute of Health and Welfare [7] |

| | | population of a given region [7]. | | |
|---|---|---|---|---|
| Healthcare Centres by Area | HCC | Healthcare centre locations across Australia calculated using the given area for a region. These healthcare centres include relevant services for children aged 0-5 [8]. Examples include Paediatric medical services, Paediatric dentistry services, and family counselling services [9]. | National Health Services Directory [9] | National Health Services Directory [9] |
| Healthcare Centres by Persons | HCC2 | Healthcare centre locations across Australia calculated using the population for a region. These healthcare centres include relevant services for children aged 0-5 [8]. Examples include Paediatric medical services, Paediatric dentistry services, and family counselling services [9]. | National Health Services Directory [9] | National Health Services Directory [9] |
| Developmentally Vulnerable on One or More Domains | DV1 | A summary indicator used to determine whether a child is developmentally vulnerability on one or more domains: physical health and well-being, social competence, | Australian Early Development Census [11] | Australian Early Development Census [11] |
| | | emotional maturity, school-based language and cognitive skills, or communication skills and general knowledge [10] | | |
| Developmentally Vulnerable on Two or More Domains | DV2 | A summary indicator used to determine whether a child is developmentally vulnerability on two or more domain: physical health and well-being, social competence, emotional maturity, school-based language and cognitive skills, or communication skills and general knowledge [10]. | Australian Early Development Census [11] | Australian Early Development Census [11] |
| Developmentally Vulnerable on the Physical Health and Well-being Domain | DVH | A child's overall physical health and well-being is measured by their readiness for the school day, gross and fine motor skills as well as their physical independence [10]. | Australian Early Development Census [11] | Australian Early Development Census [11] |
| Developmentally Vulnerable on the Emotional Maturity Domain | DVE | Social competence for children is determined through several factors. These include a child's overall social competence, their approach to learning, | Australian Early Development Census [11] | Australian Early Development Census [11] |

| | | responsibility, readiness to explore, and respect [10]. | | |
|---|---|---|---|---|
| Developmentally Vulnerable on the Social Competence Domain | DVS | Developmentally vulnerable based on relevant question and cut-off for social competence, learning and readiness, as well as responsibility and respect [10]. | Australian Early Development Census [11] | Australian Early Development Census [11] |
| Developmentally Vulnerable on the Language and Cognitive Skills (School-Based) Domain | DVL | Children's literacy and numeracy, interest in literacy is measured. This is then determined to be developmentally vulnerable based off of cut-offs and relevant questions [10]. | Australian Early Development Census [11] | Australian Early Development Census [11] |
| Developmentally Vulnerable on the Communication Skills and General Knowledge Domain | DVC | Relevant questions and cut-offs were used to determine if a child's communication skills and general knowledge was developmentally vulnerable [10]. | Australian Early Development Census [11] | Australian Early Development Census [11] |
| Developmentally Vulnerable on the Physical Readiness and Cleanliness for School Day Domain | DVR | The physical readiness and cleanliness of a child for the school day. Children were determined to be developmentally vulnerable based off of cut-offs and relevant questions [10]. | Australian Early Development Census [11] | Australian Early Development Census [11] |
| Developmentally Vulnerable on the Physical Independence Domain | DVI | Cut-offs and relevant questions were used to determine if child's physical independence was developmentally vulnerable [10]. | Australian Early Development Census [11] | Australian Early Development Census [11] |
| Developmentally Vulnerable on the Gross and Fine Motor Skills Domain | DVG | A child was determined to be developmentally vulnerable for their gross and fine motor skills based on cut-off and relevant questions [10]. | Australian Early Development Census [11] | Australian Early Development Census [11] |
| Live Low Birth Weights | LOW | Live births where the weight of the baby does not exceed 2500g [12]. Birthweight can be used as a determinant for the immediate and future health and well-being of a baby [13]. | National Perinatal Data Collection [14] | Australian Institute of Health and Welfare (National Perinatal Data Collection) [120] |
| Small for Gestational Age | SML | The birthweight is less than 2750g when a baby is born on or after 40 weeks is considered small for gestational age [118], [121]. This is a key indicator for infant health used to measure both the health and well-being of mother and is a determinant of the baby's health, survival, development and health and well-being [121]. | National Perinatal Data Collection [120] | Australian Institute of Health and Welfare (National Perinatal Data Collection) [120] |

| Mother's Smoking During Pregnancy | SMO | Women who reported smoking during pregnancy. This is a preventable risk factor that results in adverse health outcomes for mothers and their babies [122]. | Australian Child Youth and Well-being Atlas [26] | Australian Institute of Health and Welfare (National Perinatal Data Collection) [120] |
|---|---|---|---|---|
| Five Plus Antenatal Visits | ANT | Antenatal visits during pregnancy exceeds five [118]. | National Perinatal Data Collection [120] | Australian Institute of Health and Welfare (National Perinatal Data Collection) [120] |
| First Antenatal Visit is Less than Fourteen Weeks | A14 | Woman whose first antenatal visit was during their first trimester (less than fourteen weeks of gestational age). Visits during the first trimester can be associated with improved maternal health during pregnancy, and health outcomes for the child [123]. | Australian Child Youth and Well-being Atlas [26] | Australian Institute of Health and Welfare (National Perinatal Data Collection) [120] |
| First Antenatal Visit is between fourteen and nineteen weeks | A19 | First antenatal visit is after the first trimester of pregnancy and is between fourteen and nineteen weeks [123]. | Australian Child Youth and Well-being Atlas [26] | Australian Institute of Health and Welfare (National Perinatal Data Collection) [120] |
| First Antenatal Visit is over Twenty weeks | A20 | Women whose first antenatal visit is over twenty weeks [123]. | Australian Child Youth and Well-being Atlas [26] | Australian Institute of Health and Welfare (National Perinatal Data Collection) [120] |

| Birth less than thirty-seven weeks | BRTH | Babies who were born at less than thirty-seven weeks of gestations. Gestational age can be associated with a baby's overall health and generally poorer health outcomes are reported with babies that are born earlier [124]. | Australian Child Youth and Well-being Atlas [26] | Australian Institute of Health and Welfare (National Perinatal Data Collection) [120] |
|---|---|---|---|---|
| Mother's Age when giving birth is 15-19 | M15 | The age of the mother at time of birth is between 15-19. If the mother is under 20 this can be associated with adverse outcomes and increased risk of complications associated with pregnancy [125]. | Australian Child Youth and Well-being Atlas [26] | Australian Institute of Health and Welfare (National Perinatal Data Collection) [120] |
| Mother's Age when giving birth is 20-24 | M20 | The age of the mother at the time of birth is between 20-24 [125]. | Australian Child Youth and Well-being Atlas [26] | Australian Institute of Health and Welfare (National Perinatal Data Collection) [120] |
| Children who are Fully Immunised | IMM | Children who have received their immunised at 1year, 2 years, and 5 years old | Australian Child Youth and Well-being Atlas [26] | Australian Government Department of Health and Aged Care [126] |
| Fertility Rate | FRT | The number of babies born per number of woman [127] | Australian Child Youth and Well-being Atlas [26] | Australian Bureau of Statistics [127] |
| Mortality Rate per 1000 | MRT | The deaths of infants (under the age of one) per 1000 [128]. | Australian Child Youth and Well-being Atlas [26] | Australian Institute of Health and Welfare (National |

| | | | | Mortality Database) [128] |
|---|---|---|---|---|
| Mortality Rate per 100000 | MRT2 | The deaths of children aged 0 to 9 per 100000 [128]. | Australian Child Youth and Well-being Atlas [26] | Australian Institute of Health and Welfare (National Mortality Database) [128] |
| Preschool Attendance | ATT | Children that are aged four or five that have been enrolled in a preschool program [129]. | Australian Child Youth and Well-being Atlas [26] | Australian Bureau of Statistics (Preschool Education) [129] |
| Core Activity Need for Assistance | CORE | Children who have a severe or profound core activity limitation and thus need assistance in their day-to-day lives. This must limit three aspects of their daily activities (e.g., self-care, communication, or mobility) and can be long term disabilities or long-term health conditions [130]. | Australian Child Youth and Well-being Atlas [26] | Australian Bureau of Statistics (Census) [131] |
| Total Homelessness | HMLS | Homeless children living in non-convention, emergency, or short-term accommodation [132]. Australian children experiencing homelessness are amongst the most socially and economically disadvantaged. Australian Child | Australian Institute of Health and Welfare (Specialist Homelessness Services Collection) [133] | Australian Institute of Health and Welfare [7] |

| | | | Youth and Well-being Atlas [26] |
|---|---|---|---|
| Homelessness As a Result of Accommodation Reasons | ACC | The reported reason for homelessness and seeking assistance is accommodation [133]. | Australian Child Youth and Well-being Atlas [26] | Australian Institute of Health and Welfare (Specialist Homelessness Services Collection) [133] |
| Homelessness As a Result of Financial Reasons | FIN | The reported reason for homelessness and seeking assistance is financial [133]. | Australian Child Youth and Well-being Atlas [26] | Australian Institute of Health and Welfare (Specialist Homelessness Services Collection) [133] |
| Homelessness As a Result of Health Reasons | HLTH | The reported reason for homelessness and seeking assistance is health [133]. | Australian Child Youth and Well-being Atlas [26] | Australian Institute of Health and Welfare (Specialist Homelessness Services Collection) [133] |
| Homelessness As a Result of Interpersonal Relationships | IR | The reported reason for homelessness and seeking assistance is interpersonal relationships [133]. | Australian Child Youth and Well-being Atlas [26] | Australian Institute of Health and Welfare (Specialist Homelessness Services Collection) [133] |
| Couples Who Are Homeless with One Or More Children | CPLS | Clients reporting homelessness are couples with one or more children [133]. | Australian Child Youth and Well-being Atlas [26] | Australian Institute of Health and Welfare (Specialist Homelessness Services Collection) [133] |

Appendix 3. The results of the Principal Component Analysis

Table 7. The proportion of variance, cumulative variance, and eigenvalues for all vulnerability indices

| Vulnerability Index | Proportion of Variance | Cumulative Variance | Eigenvalues |
| --- | --- | --- | --- |
| VI1 | 11.19 | 43.05 | 43.05 |
| VI2 | 3.28 | 12.62 | 55.67 |
| VI3 | 2.74 | 10.52 | 66.19 |
| VI4 | 1.43 | 5.50 | 71.69 |
| VI5 | 1.10 | 4.24 | 75.94 |
| VI6 | 0.89 | 3.40 | 79.34 |
| VI7 | 0.82 | 3.17 | 82.51 |
| VI8 | 0.69 | 2.67 | 85.18 |
| VI9 | 0.63 | 2.41 | 87.60 |
| VI10 | 0.59 | 2.26 | 89.85 |
| VI11 | 0.49 | 1.90 | 91.75 |
| VI12 | 0.45 | 1.72 | 93.47 |
| VI13 | 0.35 | 1.33 | 94.81 |
| VI14 | 0.30 | 1.16 | 95.97 |
| VI15 | 0.23 | 0.88 | 96.85 |

| VI16 | 0.20 | 0.77 | 97.62 |
| VI17 | 0.17 | 0.64 | 98.26 |
| VI18 | 0.11 | 0.43 | 98.68 |
| VI19 | 0.10 | 0.37 | 99.05 |
| VI20 | 0.07 | 0.27 | 99.32 |
| VI21 | 0.05 | 0.19 | 99.51 |
| VI22 | 0.04 | 0.16 | 99.67 |
| VI23 | 0.04 | 0.14 | 99.81 |
| VI24 | 0.02 | 0.10 | 99.91 |
| VI25 | 0.02 | 0.09 | 99.99 |
| VI26 | 0.00 | 0.01 | 100.00 |

Table 8. The weights of each variable for the first five vulnerability indices (V1 to V5).

|      | VI1   | VI2  | VI3   | VI4   | VI5   |
|------|-------|------|-------|-------|-------|
| ERP  | -0.18 | 0.38 | 0.10  | -0.08 | 0.16  |
| PRIM | -0.18 | 0.37 | 0.09  | -0.06 | 0.14  |
| CC   | -0.18 | 0.39 | 0.09  | -0.01 | 0.13  |
| HCC2 | 0.05  | 0.01 | 0.12  | 0.53  | -0.36 |
| DVH  | 0.28  | 0.14 | -0.02 | -0.11 | -0.01 |

| | | | | |
|------|-------|-------|-------|-------|
| DVE | 0.26 | 0.13 | -0.05 | -0.05 | -0.06 |
| DVS | 0.26 | 0.19 | 0.00 | -0.14 | 0.06 |
| DVL | 0.27 | 0.13 | 0.01 | 0.02 | -0.08 |
| DVC | 0.23 | 0.27 | 0.09 | -0.13 | 0.01 |
| DVR | 0.27 | 0.13 | 0.04 | 0.07 | 0.01 |
| DVI | 0.24 | 0.05 | -0.08 | -0.23 | -0.13 |
| DVG | 0.25 | 0.03 | 0.01 | -0.13 | 0.19 |
| SML | 0.04 | 0.33 | 0.15 | 0.13 | 0.21 |
| SMO | 0.27 | -0.06 | 0.03 | 0.18 | 0.12 |
| ANT | -0.08 | 0.15 | -0.28 | -0.18 | -0.17 |
| A14 | -0.03 | 0.08 | -0.56 | 0.20 | 0.18 |
| A19 | -0.03 | -0.04 | 0.50 | -0.21 | 0.00 |
| A20 | 0.10 | -0.11 | 0.42 | -0.10 | -0.32 |
| BRTH | 0.23 | 0.05 | -0.03 | 0.16 | 0.01 |
| M15 | 0.25 | 0.01 | -0.02 | 0.23 | 0.01 |
| M20 | 0.26 | -0.09 | -0.03 | 0.03 | 0.08 |
| IMM | 0.08 | -0.19 | 0.05 | 0.08 | 0.53 |
| FRT | 0.21 | -0.24 | -0.02 | -0.05 | 0.08 |
| MRT2 | 0.18 | 0.28 | 0.03 | 0.16 | -0.08 |
| ATT | 0.06 | 0.09 | -0.32 | -0.37 | -0.38 |
| CORE | 0.11 | -0.21 | -0.02 | -0.40 | 0.28 |

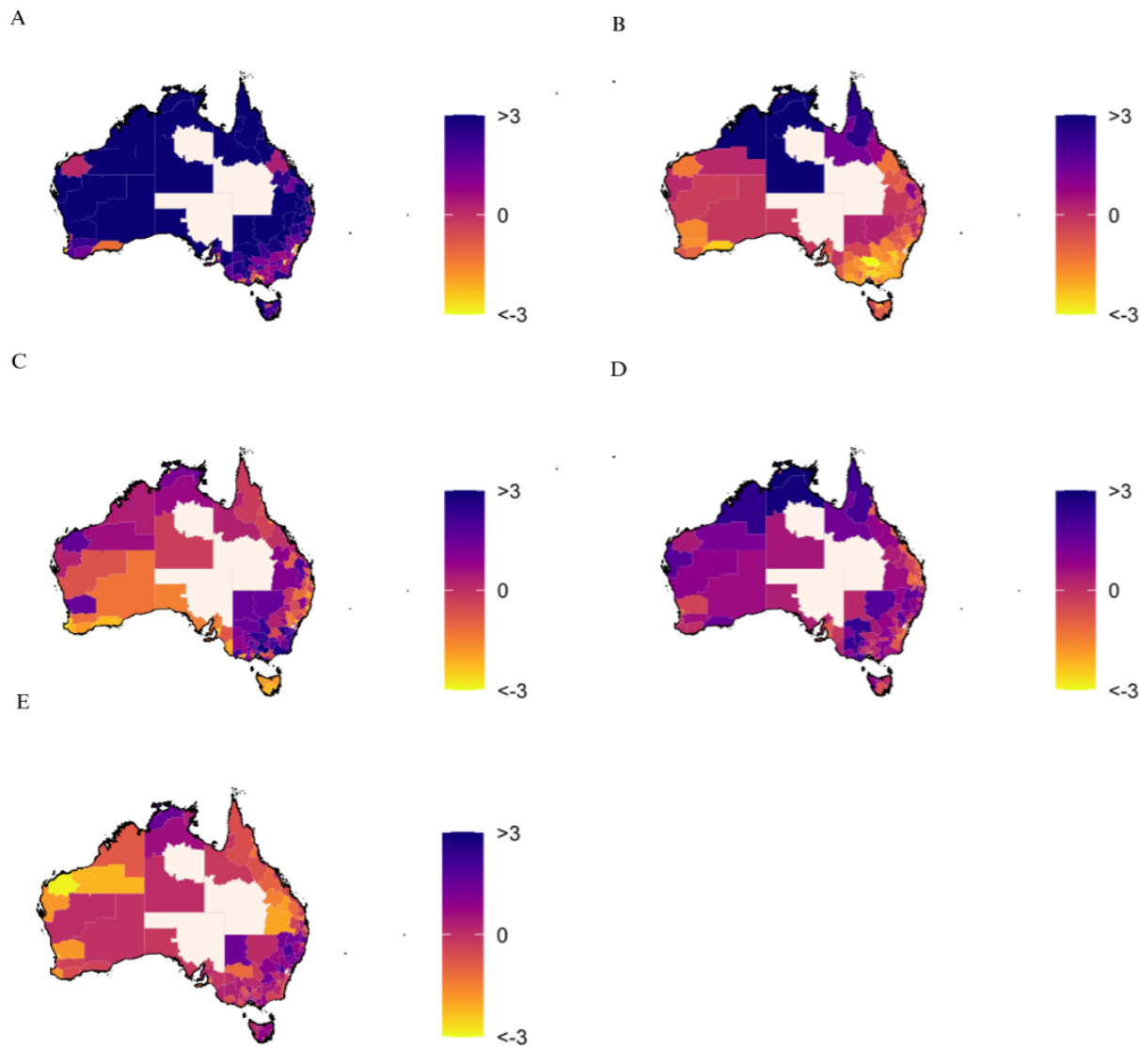

Figure 6. The alphabetical order in Figure 6 corresponds to the first five vulnerability indices in numerical order. These maps showcase each index on a consistent truncated scale. Regions with no population based on the ERP are coloured white (see Appendix 1). Purple (high positive values) indicate high levels of vulnerability and yellow (high negative values) represent levels of low vulnerability.

Appendix 4. Skewness for vulnerability indices

The level of skewness can be categorised into three different levels: ± 1 low levels of skewness (low skewness, acceptable), ± 1 to ±2 medium levels of skewness (medium skewness, acceptable), and above ± 2 level of skewness (not acceptable) (27,60). Unacceptable levels of skewness was determined to exceed ±2 (60). The first five principal components had acceptable levels of skewness (within ±2) and were used for further analysis. VI1 had a skewness of – 0.68, VI2: 0.89, VI3: -0.37, VI4: -1.42, and VI5: 0.89. The vulnerability indices have also been visualised in the following histograms.

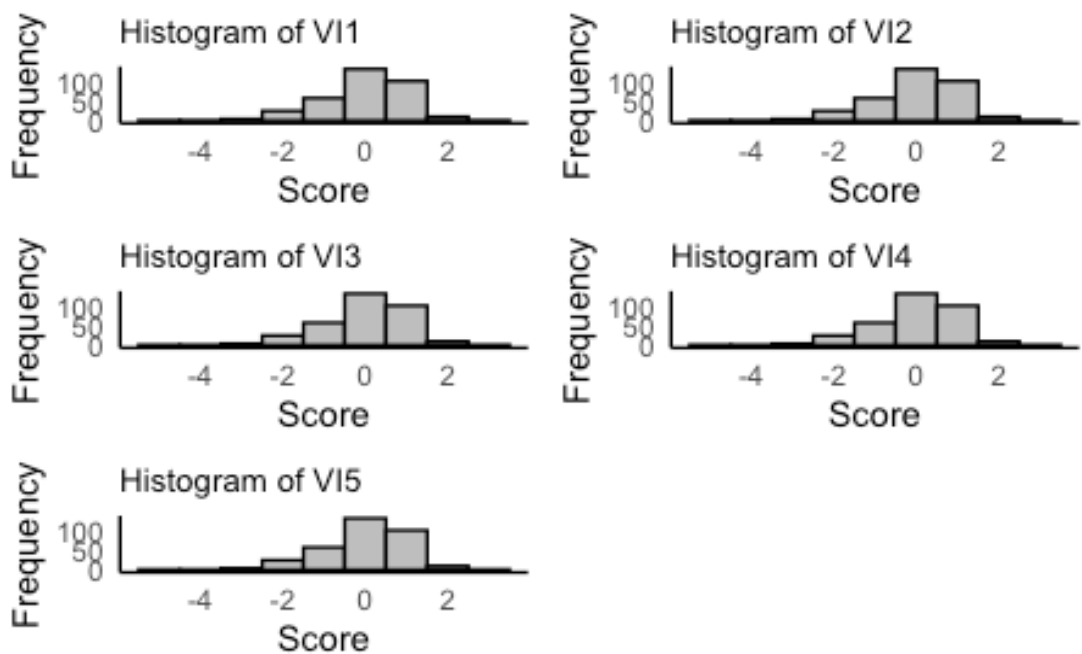

Figure 7. Histograms of the first five vulnerability indices

Appendix 5. Scree plot for the optimal number of clusters

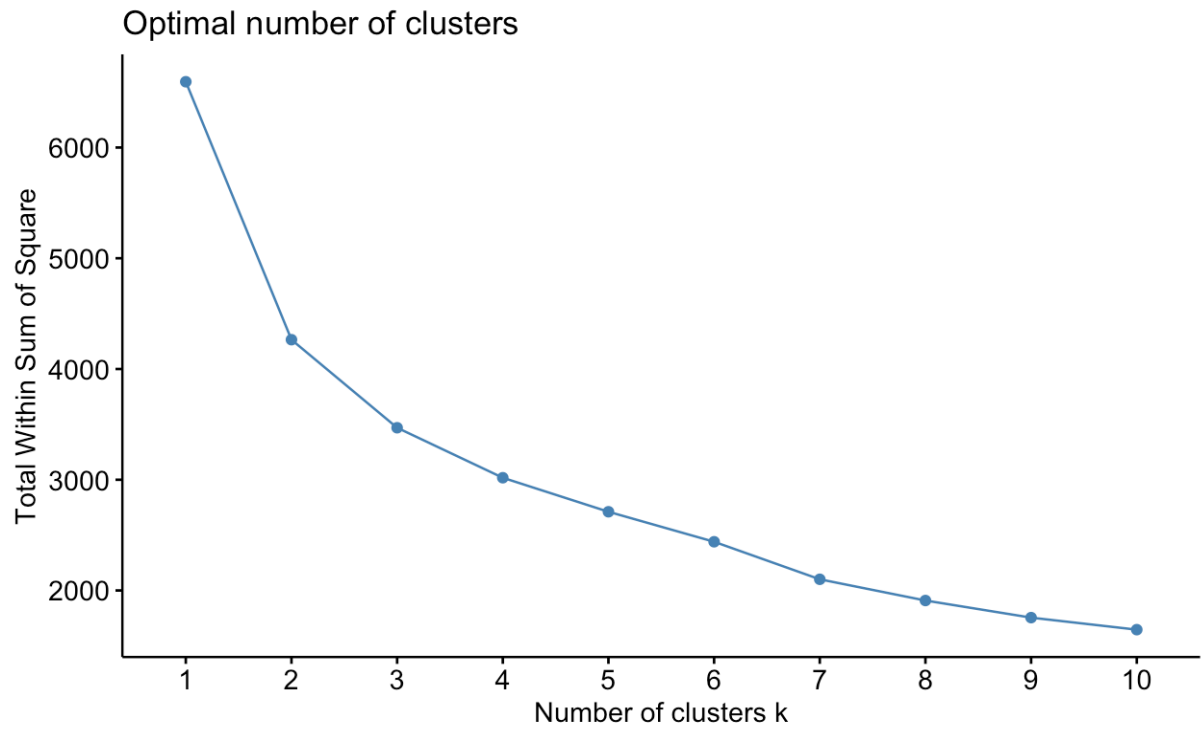